\documentclass[reqno,oneside,12pt,letterpaper]{amsart}
\usepackage{dwpreamble}
\usepackage{dwcommands}

\renewcommand{\geq}{\geqslant}

\begin{document}

\title[Lorentzian spectral zeta functions on asymptotically Minkowski spacetimes]{\large{Lorentzian spectral zeta functions \\ on asymptotically Minkowski spacetimes}\smallskip}

\author{}
\address{Institut de Mathématiques de Jussieu (UMR 7586, CNRS), Sorbonne Université -- Université de Paris, 4 pl.~Jussieu,
75252 Paris, France}
\email{nguyen-viet.dang@imj-prg.fr}
\author[]{\normalsize Nguyen Viet \textsc{Dang} \& Micha{\l} \textsc{Wrochna}}
\address{Laboratoire AGM (UMR 8088, CNRS), CY Cergy Paris Universit\'e, 2 av.~Adolphe Chauvin, 95302 Cergy-Pontoise, France}
\email{michal.wrochna@cyu.fr}

\begin{abstract} In this note, we consider perturbations of Minkowski space as well as more general spacetimes on which the wave operator $\square_g$ is essentially self-adjoint.  We review a recent result which  gives the meromorphic continuation of the Lorentzian spectral zeta function density, i.e.~of the trace density of  complex powers $\cv \mapsto (\square_g-i \varepsilon)^{-\cv}$.  In even dimension $n\geq 4$, the residue  at $\frac{n}{2}-1$  is shown to be a multiple of the scalar curvature in the limit $\varepsilon\to 0^+$. This yields a spectral action for gravity  in Lorentzian signature.
\end{abstract}

\maketitle

\section{Main result}

\subsection{Motivation} Suppose $(M,g)$ is a compact Riemannian manifold of dimension $n$, and let $\triangle_g$ be the Laplace--Beltrami operator. A classical result in analysis, dating back to  Minakshisundaram--Pleijel \cite{Minakshisundaram1949} and Seeley \cite{seeley}, states that for $\Re \cv > \n2$ the trace density of $(-\triangle_g)^{-\cv}$, defined as the on-diagonal restriction 
\beq\label{eq1}
(-\triangle_g)^{-\cv}(x,x)
\eeq
of the Schwartz kernel $(-\triangle_g)^{-\cv}(x,y)$, exists for all $x\in M$. Furthermore  \eqref{eq1} extends to a density-valued meromorphic function of the complex variable $\cv$. Its integral  over $M$  is the celebrated  \emph{spectral zeta function} of $-\triangle_g$ (or \emph{Minakshisundaram--Pleijel zeta function}),  which has attracted  widespread attention due to its  relationships with the geometry of $(M,g)$.

In fact, the   residues of  \eqref{eq1} are given by local geometric quantities: in particular if $n\geq 4$ is even, one finds
\beq \label{eq2}
\res_{\alpha=\frac{n}{2}-1} (-\triangle_g)^{-\alpha}(x,x)=\frac{{R_g(x)}}{6(4\pi)^{\n2} \Gamma(\frac{n}{2}-1) },
\eeq
where $R_g(x)$ is the scalar curvature of $(M,g)$ at $x\in M$. This identity,  often attributed to Kastler \cite{kastler} and Kalau--Walze  \cite{kalauwalze}, and  announced previously by Connes,  is  a   consequence of classical theorems in elliptic theory  (the  heat kernel based argument can be found in \cite[Thm.~1.148]{Connes2008a}; see \cite[\S1.7]{gilkey} for an approach in the spirit of Atiyah--Bott--Patodi \cite{Atiyah1975}).   Its importance in physics  stems  from the fact that the  variational equation $\delta_g R_g=0$ for $g$ is equivalent to the Einstein equations in Riemannian signature. Therefore, the l.h.s.~of \eqref{eq2} yields a \emph{spectral action} for Euclidean gravity.  Relationships of this type have also been used to justify definitions of curvature in non-commutative geometry \cite{Connes2008a,Connes2014}.

However,  it is the Einstein equations in \emph{Lorentzian} signature which have a direct physical meaning. This means that $(M,g)$ should be replaced by a Lorentzian manifold (typically not compact), but then the problem is that the corresponding Laplace--Beltrami operator $\square_g$ (or \emph{wave operator}),  is not elliptic nor bounded from below. In consequence, it is not at all clear if  $\square_g$ has a self-adjoint extension and even less clear if  the arguments from elliptic theory can be somehow replaced  (for instance, it is difficult to imagine that the heat kernel could be usefully generalized to the Lorentzian setting).

Nevertheless, it was demonstrated by Vasy \cite{vasyessential} (followed by a generalization by Nakamura--Taira  \cite{nakamurataira}) that if $(M,g)$ is well-behaved at infinity,  $\square_g$ \emph{is} essentially self-adjoint in $L^2(M,g)$.  Consequently, complex powers $(\square_g-i \varepsilon)^{-\cv}$ can be defined by functional calculus for any $\varepsilon>0$. The question  is then if this global, spectral theoretical object has  anything to do with the local geometry of $(M,g)$, in particular with the Lorentzian scalar curvature $R_g(x)$.

\subsection{Main theorem} In \cite{Dang2020} we consider Vasy's framework and provide an affirmative answer in the form of an identity largely analogous to \eqref{eq2}. Namely, we prove the following theorem. 

\begin{theorem}[{\cite[Thm.~1.1]{Dang2020}}]\label{mainthm} Assume $(M,g)$ is a globally hyperbolic, non-trapping Lorentzian scattering space of even dimension $n\geq 4$.  For all $\varepsilon>0$, the Schwartz kernel of $(\square_g-i \varepsilon)^{-\cv}$ has for $\Re\cv>\frac{n}{2}$ a well-defined on-diagonal restriction $(\square_g-i \varepsilon)^{-\cv}(x,x)$, which extends as a meromorphic function of $\cv\in\cc$ with poles at $\{\n2,$ $\n2-1$, $\n2-2$, $\dots$, ${1}\}$. Furthermore, 
\beq\label{mainresult}
\lim_{\varepsilon\to 0^+}\res_{\cv=\frac{n}{2}-1} \left(\square_g- i \varepsilon\right)^{-\cv}(x,x)=  \frac{R_g(x)}{{i}6(4\pi)^{\n2} \Gamma\big(\frac{n}{2}-1\big) }, 
\eeq
where $R_g(x)$ is the scalar curvature at $x\in M$. 
\end{theorem}

The meromorphic continuation of $\cv \mapsto \zeta_{g,\varepsilon}(\cv)(x):= \left(\square_g- i \varepsilon\right)^{-\cv}(x,x)$  is called the \emph{Lorentzian spectral zeta function density} of $(M,g)$. 

Let us briefly discuss the assumptions of Theorem \ref{mainthm}.  The class of \emph{non-trapping Lorentzian scattering spaces} introduced by Vasy \cite{vasyessential} can be thought of having asymptotically the same structure as Minkowski space at spacetime infinity $\module{x}\to +\infty$, with the extra requirement that there are no trapped null geodesics. It is worth emphasizing that this is a somewhat more general class than what one would typically call ``asymptotically Minkowski spacetime'' in that the definition refers to the bicharacteristic flow (the null geodesic flow lifted to the cotangent bundle) and to its asymptotic properties, rather than to the precise form of  the metric coefficients at infinity, see \cite{vasyessential,Dang2020}.  \emph{Global hyperbolicity} is a  standard assumption which provides a general setting for well-posedness of the Cauchy problem for $\square_g$ and is unlikely to entail significant loss of generality (in fact, it automatically follows from the non-trapping assumption for a large class of asymptotically Minkowski spacetimes, see \cite[\S{4.2}]{GWfeynman}). 

The most essential feature of these assumptions is that they allow for perturbations of Minkowski space without assuming any particular symmetries or analyticity. This means there are sufficiently many variations of the metric to derive Einstein equations from the r.h.s.~ of \eqref{mainresult}. Consequently, the l.h.s.~gives a {spectral action for gravity} in Lorentzian signature.

\subsection{Further results} Let us also briefly mention  several of our further results related to Theorem \ref{mainthm}. 

In \cite{Dang2020} we show an expansion in the spirit of the Chamseddine--Connes spectral action \cite{Connes1996,Chamseddine1997}.  Namely, for any Schwartz function $f$ with Fourier transform $\widehat{f}$ supported in $\open{0,+\infty}$ and any $N\in\nn_{\geqslant 0}$, we have for $\varepsilon>0$ the large $\scalambda>0$ expansion
\beq\label{eq:main}
 f\big((\square_g+i\varepsilon)/\scalambda^2\big)(x,x) = \sum_{j=0}^N \scalambda^{n-2j} C_j(f) \, a_j(x) + {\pazocal{O}}(\varepsilon,\scalambda^{n-2N-1}),
\eeq
where each $C_j(f)$ depends only on $j\in\nn_{\geqslant0}$, the space-time dimension $n$ and $f$, and $a_j(x)$ are directly related to the \emph{Hadamard coefficients}, in particular  
$$
a_0(x)=(4\pi)^{-\n2},  \quad a_1(x)={-}(4\pi)^{-\n2}\frac{1}{6} R_g(x), 
$$
with  $C_0(f)=i^{-1}e^{\frac{in\pi}{4}} \int_0^\infty\widehat{f}(t)t^{\frac{n}{2}-1}dt$ and   $C_1(f)=i^{-1}{e^{\frac{i(n-2)\pi}{4}} }\int_0^\infty\widehat{f}(t)t^{\frac{n}{2}-2}dt$.

Furthermore, we show that the identities \eqref{mainresult}--\eqref{eq:main} remain  valid in the case of \emph{ultrastatic spacetimes}  $(M,g)$, meaning that $M=\rr\times Y $ and  $g=dt^2-h$ for some $t$-independent complete Riemannian manifold $(Y,h)$.  In this setting essential self-adjointness is due to  Derezi\'nski--Siemssen \cite{derezinski} and the proofs are significantly simpler because the spectral theory of $-\triangle_h$ can then be used. We remark that in the related and more general case of  \emph{stationary  spacetimes} the scalar curvature can be recovered in a different spectral-theoretical way  through a Gutzwiller--Duistermaat--Guillemin trace formula due to Strohmaier--Zelditch \cite{Strohmaier2020b}.

Finally, in a further work \cite{Dang2021} we define a dynamical notion of ``residue'' which generalizes the \emph{Guillemin–Wodzicki residue density}  \cite{Guillemin1985,wodzicki} of pseudo-differential operators. More precisely, given a Schwartz kernel, our definition refers to the \emph{Pollicott--Ruelle resonances} for the dynamics of scaling towards the diagonal in $M\times M$. We apply this formalism to complex powers $(\square_g - i \varepsilon)^{-\cv}$ and we demonstrate that residues of Lorentzian spectral zeta functions $\zeta_{g,\varepsilon}(\cv)$ are dynamical residues indeed. This provides a Lorentzian version of the fact that the residue  \eqref{eq2}  can be expressed as a Guillemin--Wodzicki residue or, in physicists' terminology, a ``scaling anomaly''. 

\section{Sketch of proof} 

\subsection{From resolvent to complex powers} Let us now give a sketch of the proof of Theorem \ref{mainthm}.  Let $P=\square_g$ be the wave operator, i.e., using the notation $\module{g}=\module{\det g}$, $P$ is the differential operator
\beq\label{eq:diffop}
\bea
P&=\module{g(x)}^{-\frac{1}{2}}\partial_{x^j} \module{g(x)}^{\frac{1}{2}}g^{jk}(x)\partial_{x^k}  \\
& =  \partial_{x^j}g^{jk}(x)\partial_{x^k}+b^k(x)\partial_{x^k}
\eea
\eeq
where we sum over repeated indices, and $b^k(x)=\module{g(x)}^{-\12} g^{jk}(x)(\partial_{x^j}\module{g(x)}^{\frac{1}{2}} )$. We use the same notation $P$ for the closure of $\square_g$ acting on test functions $C_{\rm c}^\infty(M)\subset L^2(M,g)$.

\begin{figure}
 \begin{tikzpicture}[scale=1.4]
 \def\bigradius{3} \def\incangle{15}\def\littleradius{0.5}
 \draw [help lines,->] (-1.0*\bigradius, 0) -- (1.0*\bigradius,0);
 \draw [help lines,->] (0, -0.5*\bigradius) -- (0, 1.0*\bigradius);
 \begin{scope}[shift={(0,2*\littleradius)}]
  \node at (-0.9,0.42) {$\gamma_{\varepsilon}$};
\path[draw,line width=0.8pt,decoration={ markings,
      mark=at position 0.15 with {\arrow{latex}}, 
      mark=at position 0.53 with {\arrow{latex}},
      mark=at position 0.85 with {\arrow{latex}}},postaction=decorate] (-\bigradius,{\bigradius*tan(\incangle)})   -- (-\incangle:-\littleradius) arc (180-\incangle:360+\incangle:\littleradius)   -- (\bigradius,{\bigradius*tan(\incangle)});
      \path[draw,line width=0.2pt,postaction=decorate,<->] (0,0)   -- (330:\littleradius);
\end{scope}
\path[draw,line width=0.2pt,postaction=decorate,dashed] (-\bigradius,2*\littleradius)   -- (\bigradius,2*\littleradius);
\path[draw,line width=0.2pt,postaction=decorate,<->] (2,0)   -- (2,2*\littleradius);
\node at (2.15,\littleradius){\scalebox{0.8}{$\varepsilon$}};
\node at (0.2,1.45*\littleradius){\scalebox{0.8}{$\frac{\varepsilon}{2}$}};
 \node at (2.9,-0.2){$\scriptstyle \Re z$};
 \node at (0.35,2.8) {$\scriptstyle i \Im z$};
 \end{tikzpicture}
 \caption{\label{fig:contour}The contour $\gamma_\varepsilon$ used to express $(P-i\varepsilon)^{-\cv}$ as an integral of the resolvent.} 
 \end{figure}
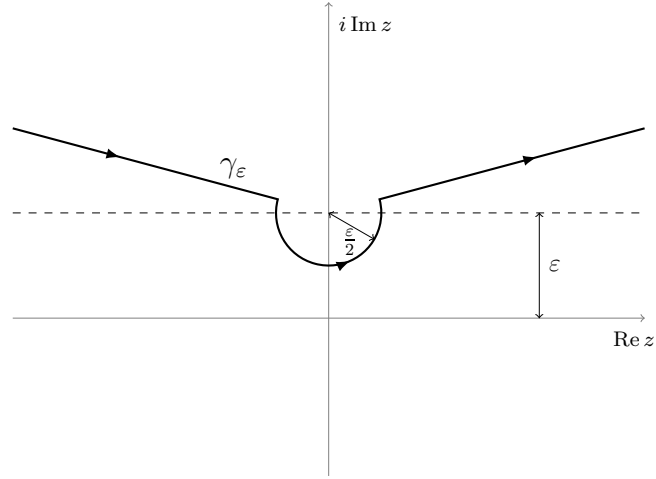

For $\varepsilon>0$ and $\Re \cv > 0$, the  power $(P-i \varepsilon)^{-\cv}$ can be expressed as a contour integral of the form
\beq\label{eq3}
(P-i\varepsilon)^{-\cv}=\frac{1}{2\pi i}\int_{\gamma_\varepsilon} (z-i\varepsilon)^{-\cv} (P-z)^{-1} dz,
\eeq
convergent in the strong operator topology (see e.g.~ \cite[App.~B]{Dang2020}).  The contour of integration  $\gamma_\varepsilon$ is represented in Figure \ref{fig:contour} and can be written as $\gamma_\varepsilon= \tilde\gamma_\varepsilon+i\varepsilon$, where
      \beq
      \tilde\gamma_{\varepsilon} = e^{i(\pi-\theta)}\opencl{-\infty,\textstyle\frac{\varepsilon}{2}}\cup \{\textstyle\frac{\varepsilon}{2} e^{i\omega}\, | \, \pi-\theta<\omega<\theta\}\cup  e^{i\theta}\clopen{\textstyle\frac{\varepsilon}{2},+\infty}
      \eeq
       goes from $\Re z\ll 0$ to  $\Re z\gg 0$ in the upper half-plane (for some fixed  $\theta\in\open{0,\pid}$).
       
 The strategy is then to construct a sufficiently explicit parametrix for the resolvent $(P-z)^{-1}$.  When estimating error terms, a significant difficulty is the necessity to control what happens \emph{uniformly in $z$}, with an appropriate decay rate along the infinite contour $\gamma_\varepsilon$. We remark that retarded and advanced propagators for $P-z$ are \emph{not} expected to have this kind of decay,  so in practice it is not possible to use various techniques from hyperbolic PDEs related to solving a retarded or advanced problem or a Cauchy problem for $P-z$.

\subsection{Uniform Hadamard parametrix}\label{ss:uhp} In contrast to the heat kernel, the \emph{Hada\-mard parametrix} for the Laplace--Beltrami operator generalizes well to the Lorentzian case. Furthermore, it is known to have  similar local geometric content.  So, the main question is whether the Hadamard parametrix approximates the resolvent $(P-z)^{-1}$ in a reasonable sense, uniformly in $z$ along the contour $\gamma_\varepsilon$.
 
Before answering this question, let us recall the construction of the Hadamard parametrix for $P-z$. 

 As expected  from explicit formulae on Minkowski space and from the theory of Fourier integral operators, there are actually four different Hadamard parametrices with different singularities. In the case of the resolvent $(P-z)^{-1}$ with $\Im z > 0$, one expects that the \emph{Feynman} Hadamard parametrix is the correct choice, see e.g.~\cite{Lewandowski2020} for the general definition. Here we use a construction directly adapted from earlier works in the Riemannian or Lorentzian time-independent case \cite{HormanderIII,soggeHangzhou,StevenZelditch2012,Zelditch2017} (cf.~\cite{Baer2020} for a unified treatment of even and odd dimensions), supplemented by new estimates that are uniform in $z$ (cf.~\cite{Sogge1988,Ferreira2014,Bourgain2015} for uniform estimates in the Riemannian case). As expected, their proof  is significantly complicated by  light-cone singularities not present in the Riemannian analogue of the problem.

\subsubsection*{Step 1} Let $\eta=dx_0^2-(dx_1^2+\cdots+dx_{n-1}^2)$ be the  Minkowski metric on $\rr^n$, and consider the corresponding quadratic form $$
\vert\xi\vert_\eta^2 =  -\xi_0^2+\textstyle\sum_{i=1}^{n-1}\xi_i^2,
$$
defined for convenience with a {minus} sign. For $\cv\in \cc$ and $\Im z >0$,  the distribution $\left(\vert\xi\vert_\eta^2-z \right)^{-\cv}$ is well-defined by pull-back from $\rr$. More generally, for $\Im z \geqslant 0$, the limit $\left(\vert\xi\vert_\eta^2-z -i0\right)^{-\cv}=\lim_{\varepsilon\to 0^+}\left(\vert\xi\vert_\eta^2-z -i\varepsilon\right)^{-\cv}$ from the upper half-plane is well defined as a distribution on $\rr^n\setminus \{0\}$. If $z\neq 0$ it can be  extended to a  family of homogeneous\footnote{Homogeneity refers here to  rescaling simultaneously the $\xi$ variables by $\lambda>0$  and the complex number $z$ by $\lambda^2$.} distributions on $\rr^n$, holomorphic in $\cv\in \cc$.  We introduce special notation for its appropriately normalized inverse Fourier transform,
\beq\label{eq:defFsz}
\Fs{z}\defeq\frac{\Gamma(\cv+1)}{(2\pi)^{n}} \int e^{i\left\langle x,\xi \right\rangle}\left(\vert\xi\vert_\eta^2-i0-z \right)^{-\cv-1}d^{n}\xi.
\eeq

\subsubsection*{Step 2}   Next, one pull-backs the distribution $\Fse{z}$ to a neighborhood $\cU$ of the diagonal $\diag\subset M\times M$  using the exponential map. In view of the  $O(1,n-1)_+^\uparrow$-invariance of  $\Fse{z}$   there is a canonical way to define this pull-back (see \cite[\S5.1]{Dang2020}), denoted in the sequel by $\Fe{z}$.  
  The Hadamard parametrix (or rather its Schwartz kernel) is  constructed in normal charts    using the family $\Fe{z}$. Namely, for fixed
   $\varvarm\in M$,  one  expresses the distribution $x\mapsto \mathbf{F}_\cv(z,\varvarm,x)$
in normal coordinates centered at $\varvarm$, defined on some $U\subset T_{x_0}{M}$.  By abuse of notation we continue to write  $\Fe{z}$ instead of $ \mathbf{F}_\cv(z,\varvarm,\exp_\varvarm(\cdot))\in \pazocal{D}^\prime(U)$. One then defines for large $N$ a parametrix $H_N(z,.)$ by setting
\begin{eqnarray}
{H_N(z,.)=\sum_{k=0}^N u_k \Fe[k]{z} \in \pazocal{D}^\prime(U) },
\end{eqnarray}
where $(u_k)_{k=0}^\infty$  is   a sequence of functions 
in $C^\infty(U)$ that solves the hierarchy of transport equations
\begin{eqnarray}
{2k u_k+ b^i(x)\eta_{ij}x^j u_k+2 x^i\partial_{x^i} u_k+2Pu_{k-1}=0}
\end{eqnarray}
with initial condition $u_0(0)=1$ (by convention, $u_{k-1}=0$ for  $k=0$,  we sum over repeated indices, and we recall that $b^i(x)$ is defined in \eqref{eq:diffop}).  The transport equations imply that $H_N(z,.)$ solves
\begin{equation}
\left(P-z\right)H_N(z,.)={\module{g}^{-\frac{1}{2}}}\delta_0+(Pu_N)\mathbf{F}_N,
\end{equation}
on  $U$, where $(Pu_N)\mathbf{F}_N$ is interpreted as an error term. 

\subsubsection*{Step 3} In the final step one takes into account the dependence on $x_0$ to obtain a parametrix on the neighborhood $\cU$ of the diagonal. Here we make this step implicitly by sticking to the same  notation $\Fe{z}$ for the corresponding distribution on $\cU$. Finally, one uses a  cutoff function $\chi\in \cf(M^2)$ supported in $\pazocal{U}$ (with $\chi=1$ near the diagonal) to extend the definition of $H_N(z,.)$ to $M^2$:
$$
\boxed{H_N(z,.)=\sum_{k=0}^N \chi u_k \Fe[k]{z}\in \pazocal{D}^\prime({M\times M}).}
$$
The Hadamard parametrix extended to $M^2$  satisfies 
\beq\label{eq:PzHN}
\left(P-z\right)  H_N(z,.)={ \module{g}^{-\frac{1}{2}}}\delta_{\diag}+(Pu_N)\mathbf{F}_N(z,.)\chi+r_N(z,.),
\eeq
where 
 $\module{g}^{-\frac{1}{2}}\delta_\Delta(x_1,x_2)$ is the Schwartz kernel of the identity map and $r_N(z,.)\in  \pazocal{D}^\prime({M\times M})$ is an error term supported in a punctured neighborhood of $\Delta$ which is due to the presence of the cutoff $\chi$.
 \medskip
 
 In order to conclude a relationship between the resolvent $(P-z)^{-1}$ and the Hada\-mard parametrix $H_N(z,.)$, the natural next step is to apply $(P-z)^{-1}$ to both sides of \eqref{eq:PzHN}. The objective is then to show that the composition of $(P-z)^{-1}$  with the two error terms on the r.h.s.~exists, decreases in $z$ in a suitable sense along the contour $\gamma_\varepsilon$, and is sufficiently  regular (so that its on-diagonal restriction always exists and the corresponding integral on $\gamma_\varepsilon$ is holomorphic in $\cv$). 

It turns out that by choosing $N$ sufficiently high we can make  $(Pu_N)\mathbf{F}_N(z,.)\chi$ decaying in $z$ and of arbitrarily high Hölder regularity. The proof is quite technical as it uses oscillatory integral representations, but most of the analysis is carried   out on the level of  the explicit model family $\Fse{z}$. In combination with regularity properties of $(P-z)^{-1}$ obtained as a corollary of Vasy's proof of essential self-adjointness, this yields an easily controllable  error term.

On the other hand, the error term $r_N(z,.)$ (although it can be arranged to be supported away from the diagonal) is always \emph{singular} regardless of the choice of $N$. This stands in sharp contrast with analogous constructions in the Riemannian case and is the most significant obstacle in the proof: a priori it is  not even clear if  the composition $(P-z)^{-1} r_N$ makes sense. 

A way out is possible thanks to a remarkable property shared by the Feynman Hadamard parametrix and $(P-z)^{-1}$ when $\Im z >0$. Their Schwartz kernels are singular, but in a special way which allows operator composition nevertheless, and which implies that the compositions have singularities of the same type.  Microlocally, they behave as the \emph{Feynman propagator} on Minkowski space, i.e.~the Fourier multiplier by $(-\xi_0^2 + \xi_1^2 + \cdots+ \xi_{n-1}^2  - i0)^{-1}$.  This condition  can be formulated in terms of an operatorial Sobolev wavefront set $\wf^{'(s)}\big((P-z)^{-1}\big)$ for large $s\in \rr$: by definition, a pair of points $(q_1,q_2)\in (T^*M\setminus \zero)^{\times 2}$ does \emph{not} belong to $\wf^{'(s)}\big((P-z)^{-1}\big)$ if there exists pseudo-differential operators $B_1,B_2\in \Psi^0(M)$, elliptic at respectively $q_1,q_2$, such that 
$$
B_1 (P-z)^{-1} B_2^* : H^{m}_\c(M)\to H^{m+s}_\loc (M)
$$
is bounded  for all $m\in \rr$. A uniform version particularly well adapted to our needs can be defined by requiring that the operator semi-norms are $O(\bra z \ket^{-\12})$ in $z$ along the contour $\gamma_\varepsilon$.

For fixed $z$, it is relatively easy to find the wavefront set of $H_N(z)$, and one could try various existing techniques to estimate the wavefront set of $(P-z)^{-1}$. However, estimates on the  \emph{uniform} wavefront set are needed to control the contributions of the error terms after integration. The uniform wavefront set of $r_N$ is obtained from a  detailed Hölder regularity analysis of  oscillatory integral representations with the help of dyadic decompositions. The uniform wavefront set of $(P-z)^{-1}$ is estimated in several steps outlined in the next paragraphs,  with a central role played by  \emph{microlocal propagation estimates} including \emph{radial estimates}. Uniform regularity of the composition    $(P-z)^{-1} r_N$ is then deduced from the   two uniform wavefront sets and the property that $r_N$ is supported away from the diagonal in $M\times M$. 

\subsection{Uniform microlocal resolvent estimates}  In the estimation of the uniform wavefront set of $(P-z)^{-1}$, the first step is to construct a parametrix  $G_z=G_z^+ + G_z^-$ for $(P-z)^{-1}$, which  consists of two terms $G_z^\pm$ that correspond each to solving an evolution problem of \emph{first order in time}. This parametrix is used as reference operator with more easily computable wavefront set.  

The construction relies on an approximate factorization of $P-z$. Namely, we show that after  a suitable  coordinate  change $\varphi$ and  a conformal transformation by some smooth factor $c>0$ (this step uses global hyperbolicity), $P-z$ can be written in  the form
$$
\bea
-c^2  (\varphi^* (P-z)) &= (D_t - A(t,z))(D_t+B(t,z)) +   R(t,z) \\
&=(D_t +\tilde B(t,z))(D_t - \tilde A(t,z)) +   \tilde R(t,z),
\eea
$$
where $A(t,z),B(t,z),\tilde A(t,z),\tilde B(t,z)\in \Psi^{1}(M)$ are smooth (in $t$) families of pseudo-differential operators which are \emph{elliptic with parameter}  in the sense of Shubin's  parame\-ter-dependent calculus  \cite{shubin}, with positive  principal symbols,  and $R(t,z),\tilde R(t,z)$ are smooth families of operators with arbitrarily good regularity properties,  uniformly in $z$. The operators $G_z^\mp$ are defined through an expression which uses  the retarded problem  of $D_t - \tilde A(t,z)$, resp.~advanced problem for $D_t+B(t,z)$ (these are the only two that are well-behaved for large $\Im z >0$). As such, their uniform wavefront sets can be estimated by arguments closely related to Egorov's theorem.

The problem is then how to demonstrate that $G_z=G_z^+ + G_z^-$  and $(P-z)^{-1}$ have the same uniform wavefront set. Since the wavefront sets of $G_z^+$ and $G_z^-$ are disjoint (they propagate singularities in the two different components $\Sigma^\mp$ of the characteristic set  of  $P$), it actually suffices to estimate the wavefront set of $(P-z)^{-1}- G_z^\pm$. The key ingredient are \emph{microlocal propagation estimates}, which can be applied  if we  have some microlocal regularity  of   $(P-z)^{-1}- G_z^\pm$ to start with.

It turns out that there is indeed a significant  property shared by $(P-z)^{-1}$ and $G_z^\pm$.    Let us first explain  it in the case of the resolvent $(P-z)^{-1}$.  Its mapping properties are best understood  in the framework of anisotropic \emph{scattering Sobolev spaces} $\Hsc{{s},{\ell}}$: these spaces generalize the weighted Sobolev spaces  $( 1+ \module{x}^2 )^{-\ell/2} H^{s}(\rr^n)$ in a way that allows the weight orders $\ell$ to vary in phase space (more specifically, on Melrose's \emph{scattering bundle} ${}^{\rm sc}T^*M$ \cite{melrosered}, which in our context provides the natural framework for microlocal analysis on the compactification of $M$). The key ingredient in Vasy's proof of essential self-adjointness  is a  \emph{Fredholm estimate} of the form 
\beq\label{est1}
\|  u\|_{s,\ell} + (\Im z)  \|  u\|_{s-\12,\ell+\12}  \leqslant C (\|  (P-z) u \|_{s-1,\ell+1} + \err{u}), 
\eeq
uniformly for $z\in \gamma_\varepsilon$ (with $\err{u}$ representing a negligible error term). Here, $\ell$ is chosen monotone along the bicharacteristic flow, in such a way that $\ell > -\12$ at \emph{sources at infinity} (from which bicharacteristics are assumed to originate), and $\ell<-\12$ at \emph{sinks at infinity} (to which bicharacteristics tend). One can think of  this condition as imposing boundary conditions at infinity: solving $(P-z)u=f$ in the corresponding spaces is then a  \emph{Feynman problem} \cite{GHV,Gerard2020,vasyessential,Taira2020a}.    The  estimate \eqref{est1} is responsible for the fact that if  $f=(P-z)u$ is compactly supported  then it  \emph{decays  at a rate faster than the threshold value $-\12$} microlocally at the sources. This statement can be improved in various ways, and $(P-z)^{-1} f$ has of course even better decay properties. The key point is     that  within $\Sigma^\mp$, $G^\pm f$ is decaying at the same source as $(P-z)^{-1}f$.
Therefore, $((P-z)^{-1}- G_z^\pm)f$  decays at  the sources microlocally in the respective component, and this property enables the use of \emph{radial estimates} in  Melrose's scattering calculus $\Psi_\sc(M)$ \cite{melrosered,vasygrenoble,vasyessential}, to get high regularity of $A((P-z)^{-1}- G_z^\pm)f$  if $A\in \Psi^{0,0}_\sc(M)$ is microsupported near sources in the respective component (incidentally, these are the same estimates which are used to prove \eqref{est1}). Then, propagation of singularities and elementary manipulations with operatorial wavefront sets are used to deduce that $(P-z)^{-1}- G_z^\pm$ is  (everywhere) smoothing. Crucially, in each step of this proof, the uniformity in $z$ is under control.

This proves the desired estimate on the uniform  wavefront set of $(P-z)^{-1}$, and  as explained in \sec{ss:uhp},  concludes the proof that $(P-z)^{-1}$ equals the Feynman Hadamard parametrix $H_N(z)$ modulo inessential terms.

\subsection{Extraction of the scalar curvature}   From that point on we can effectively replace the resolvent $(P-z)^{-1}$ with the Hadamard parametrix $H_N(z)$.   In fact, if we  integrate   $(z-i\varepsilon)^{-\cv} H_N(z)$ over the contour $\gamma_\varepsilon$ instead of $(z-i\varepsilon)^{-\cv}(P-z)^{-1}$, the result will   differ  from $(P-i\varepsilon)^{-\cv}$ merely by a term whose trace density is holomorphic in $\cv$. 

The integral turns out to be of the same form as the Hadamard expansion. More precisely, $(P-i\varepsilon)^{-\cv}$ equals 
\beq\label{rrr}
 \sum_{\varm=0}^N \chi u_\varm\frac{(-1)^\varm\Gamma(-\cv+1)}{\Gamma(-\cv-\varm+1)\Gamma(\cv+\varm)} \mathbf{F}_{\varm+\cv-1}(-i\varepsilon,.).
 \eeq
plus the irrelevant error term.   The meromorphic properties of the on-diagonal restriction of \eqref{rrr} can be deduced from an analysis on $\rr^n$ thanks to the identity $\mathbf{F}_{\cv}(z,x,x)=\Fso{z}$, valid for every $x\in M$.

A toy example illustrating what happens in Euclidean signature is provided by  the integral 
$$
\int_{\mathbb{R}^n}(\Vert\xi\Vert^2-z)^{-\cv}d^n\xi=\frac{1}{\Gamma(\cv)} \int_0^\infty \left( \int_{\mathbb{R}^n} e^{-t(\Vert\xi\Vert^2-z)}d^n\xi\right) t^{\cv-1} dt,
$$
assuming for simplicity $z<0$ for the moment.
It has the same poles as
$$ \bea 
&\frac{1}{\Gamma(\cv)}  \int_0^1 \left( \int_{\mathbb{R}^n} e^{-t(\Vert\xi\Vert^2-z)}d^n\xi\right) t^{\cv-1} dt \\ &=\frac{(2\pi)^n}{\Gamma(\cv)(4\pi)^{\frac{n}{2}}}
 \sum_{k=0}^\infty \frac{z^k}{k!} \int_0^1  t^{\cv-\frac{n}{2}+k-1} dt=\frac{\pi^{\frac{n}{2}}}{\Gamma(\cv)}
 \sum_{k=0}^\infty \frac{z^k}{k!(\cv-\frac{n}{2}+k)}. 
\eea $$
In consequence, we see that the residue at $\cv=k$, $k\in \{1,\dots,\frac{n}{2}-1\}$
is
\begin{equation}
\label{eq:yot}
\Res_{\cv=k}\int_{\mathbb{R}^n}(\Vert\xi\Vert^2-z)^{-\cv}d^n\xi=\frac{z^{\frac{n}{2}-k}\pi^{\frac{n}{2}}}{(\frac{n}{2}-k)! \Gamma(k)}.
\end{equation}

In our problem, we  need  to deal with integrals involving the Minkowski quadratic form rather than the Euclidean one. To that end we consider the  complex valued $n$-form
$$\omega_{\cv}=\bigg(\sum_{i=1}^n \xi_i^2-z\bigg)^{-\cv}d\xi_1\wedge\dots\wedge d\xi_n\in \Omega^{n,0}$$
for $z$ in the upper half-plane. We show that it is closed, and that  Stokes' theorem can be applied to deform the signature from Euclidean to Lorentzian in integrated expressions, which eventually yields 
$$
\bea
\Res_{\cv=k}\int_{\mathbb{R}^n} \bigg(-\xi_1^2+  \sum_{i=2}^n \xi_i^2-z-i0\bigg)^{-\cv}d^n\xi&=i \Res_{\cv=k}\int_{\mathbb{R}^n}\bigg(\sum_{i=1}^n \xi_i^2-z\bigg)^{-\cv}d^n \xi, 
\eea
$$
where the r.h.s.~is computed using \eqref{eq:yot}.

By taking into account the $\Gamma$ function factors in \eqref{rrr} we get the location of the poles of $(P-i\varepsilon)^{-\cv}$, and the remaining ingredient in the computation of the residues are the on-diagonal restrictions $u_k(x,x)$ of the  coefficients $u_k(x,y)$. These coefficients can be found for instance by observing that the transport equations for $u_k$ are analogous to transport equations in the Riemannian setting, so they are given by analogous expressions in terms of the metric $g$ and its derivatives, with obvious sign changes to account for the switch of signature \cite{Moretti1999,Bytsenko2003}. The Riemannian transport equations are in turn directly related to transport equations for the more familiar  heat kernel coefficients.

We are particularly interested in the residue at $\cv=\frac{n}{2}-1$ which comes from the coefficient $u_1(x,x)$,  and this coefficient can also be found by an inspection of the first two transport equations, directly in Lorentzian signature.  In normal coordinates (also denoted   by $x$) centered around an arbitrary point $x_0\in M$ (so  $x_0$ is $x=0$ in normal coordinates),  we have the identity
$$
\boxed{P=\partial_{x^k} g^{kj}(x)\partial_{x^j}+g^{jk}(x)(\partial_{x^j}\log \module{g(x)}^{\frac{1}{2}} ) \partial_{x^k} }.
$$
This can be used to  express the transport equations in a more convenient form and one finds after a short computation that they imply
 $$
 u_1(0)=-Pu_0{(0)}=- P(\module{g({0})}^\frac{1}{4} \module{g(x)}^{-\frac{1}{4}})|_{x=0}. $$ 
In normal coordinates, $\vert g(0)\vert^\frac{1}{4}=1$ and
$$g_{ij}(x)=\eta_{ij}+\frac{1}{3}R_{ikjl}x^kx^l +\pazocal{O}(\vert x\vert^3), \quad 
\module{g(x)}^{-\frac{1}{4}}=1+\frac{1}{12}\mathbf{Ric}_{kl}(0)x^kx^l+\pazocal{O}(\vert x\vert^3),
$$
where $\mathbf{Ric}_{kl}$ is the {Ricci} tensor. This implies that
$$-P \module{g(x)}^{-\frac{1}{4}}=-\frac{1}{6} { g^{kl}\mathbf{Ric}_{kl}(0)}+\pazocal{O}(\vert x\vert), 
$$ 
where 
$g^{kl}\mathbf{Ric}_{kl}=R_g(0)$ is the scalar curvature at $x_0$. Since $x_0$ was arbitrary, we conclude  $u_1(x,x)=-\frac{1}{6}R_g(x)$.

\medskip

\subsubsection*{Acknowledgments} The authors gratefully acknowledge support from the grant ANR-20-CE40-0018. 

\bibliographystyle{abbrv}
\bibliography{complexpowers}

\renewcommand{\geq}{\geqslant}\end{document}